\documentclass[11pt]{amsart}
\usepackage{amssymb, color}
\usepackage{amsmath}
\usepackage{amsthm}
 \font\tenmsb=msbm10 at 12pt \font\sevenmsb=msbm7 at 8pt \font\fivemsb=msbm5 at
6pt
\newfam\msbfam
\textfont\msbfam=\tenmsb \scriptfont\msbfam=\sevenmsb \scriptscriptfont\msbfam=\fivemsb
\def\Bbb#1{{\tenmsb\fam\msbfam#1}}

\def\N{\Bbb N}

\def\S{\Bbb S}

\begin{document}
\def \theequation{\arabic{section}.\arabic{equation}}
\newcommand{\reset}{\setcounter{equation}{0}}

\newcommand{\beq}{\begin{equation}}
\newcommand{\noi}{\noindent}
\newcommand{\eeq}{\end{equation}}
\newcommand{\dis}{\displaystyle}
\newcommand{\mint}{-\!\!\!\!\!\!\int}

\def \bx{\hspace{2.5mm}\rule{2.5mm}{2.5mm}} \def \vs{\vspace*{0.2cm}} \def
\hs{\hspace*{0.6cm}}
\def \ds{\displaystyle}
\def \p{\partial}
\def \O{\Omega}
\def \b{\beta}
\def \m{\mu}
\def \l{\lambda}
\def \le{\lambda^*}
\def \ul{u_\lambda}
\def \D{\Delta}
\def \d{\delta}
\def \s{\sigma}
\def \e{\varepsilon}
\def \a{\alpha}
\def \g{\gamma}
\def \R{\mathbb{R}}
\def \S{\mathbb{S}}
\def\qed{%
\mbox{ }%
\nolinebreak%
\hfill%
\rule{2mm} {2mm}%
\medbreak%
\par%
}
\newtheorem{thm}{Theorem}[section]
\newtheorem{lem}[thm]{Lemma}
\newtheorem{cor}[thm]{Corollary}
\newtheorem{prop}[thm]{Proposition}
\theoremstyle{definition}
\newtheorem{defn}[thm]{Definition}
\newtheorem{rem}[thm]{Remark}
\newenvironment{thmskip}{\begin{thm}\hfill}{\end{thm}}

\newcommand{\Om}{\Omega}

\def\p{\partial}
\def \l{\lambda}
\def \pr {\noindent {\it Proof:} }
\def \rmk {\noindent {\it Remark} }
\def \esp {\hspace{4mm}}
\def \dsp {\hspace{2mm}}
\def \ssp {\hspace{1mm}}

\title{Qualitative properties of H\'enon type equations with exponential nonlinearity}
\author{Zongming Guo}
\address{Department of Mathematics, Henan Normal University, Xinxiang, 453007, China}
\email{gzm@htu.cn}

\author{Xia Huang}
\address{School of Mathematical Sciences, East China Normal
University, Shanghai, 200241, China}
\email{xhuang@cpde.ecnu.edu.cn}
\author{Dong Ye}
\address{School of Mathematical Sciences, East China Normal
University, Shanghai, 200241, China}
 \email{dye@math.ecnu.edu.cn}
\author{Feng Zhou}
\address{School of Mathematical Sciences, East China Normal
University, Shanghai, 200241, China}
\email{fzhou@math.ecnu.edu.cn}
\date{}
\begin{abstract}
We are interested in the qualitative properties of solutions of
the H\'enon type equations with exponential nonlinearity. First, we classify the stable at infinity solutions of $\Delta u +|x|^\alpha e^u=0$ in $\R^N$, which gives a complete answer to the problem considered in \cite{WY}. Secondly, existence and precise asymptotic behaviors of entire radial solutions to $\Delta^2 u=|x|^{\alpha} e^u$ are obtained. Then we classify the stable and stable at infinity radial solutions to $\Delta^2 u=|x|^{\alpha} e^u$ in any dimension.
\end{abstract}

\maketitle

\section{introduction}
The main objective here is to understand stable at infinity solutions of the second and fourth order H\'enon type equations with exponential nonlinearity, namely
\begin{equation}
\label{new1.1}
\Delta u+|x|^\alpha e^u=0 \;\;\; \mbox{in $\mathbb{R}^N$},
\end{equation}
and
\begin{equation}
\label{1.1-}
\Delta^2 u=|x|^{\alpha} e^u \;\;\; \mbox{in $\mathbb{R}^N$}.
\end{equation}
where $N \geq 2$ and $\alpha > -2$.

\medskip
Let's begin with \eqref{new1.1}. By a solution to \eqref{new1.1}, we mean that $|x|^\alpha e^u \in L_{loc}^1 (\R^N)$, $\nabla u \in L^2_{loc}(\R^N)$ and
\begin{equation*}
\int_{\R^N} \big(\nabla u \cdot \nabla \psi-|x|^\alpha e^u \psi\big) dx=0, \;\;\; \forall\; \psi \in C_0^1 (\R^N).
\end{equation*}

A solution $u$ of \eqref{new1.1} is said to be {\it stable} on a domain $\Omega \subset \R^N$ if
\begin{equation*}
Q_u (\varphi):=\int_\Omega \Big(|\nabla \varphi|^2- |x|^\alpha e^u \varphi^2 \Big) \geq 0, \;\;\; \forall\; \varphi \in C^1_0 (\Omega).
\end{equation*}
The Morse index of a solution $u$ on $\Om$, ${\rm ind}_{\Om}(u)$ is defined as the maximal dimension of all subspaces $X$ in $C_0^1(\Omega)$ such that $Q_u (\psi)<0$ for any $\psi \in X \backslash \{0\}$. Readily $u$ is stable on $\Om$ if and only if ${\rm ind}_{\Om}(u) = 0$. We say that $u$ is stable at infinity, when $u$ is stable out of a compact set. It's well known that any finite Morse index solution over $\R^N$ is stable at infinity.

\medskip
For $\alpha=0$, Farina \cite{Fa1} showed that $\Delta u+e^u=0$ has no stable classical solution in $\R^N$ if $2 \leq N \leq 9$.
He proved also that any classical solution with finite Morse index in $\R^2$ verifies $e^u \in L^1 (\R^2)$, so
must be a sphere solution following Chen \& Li \cite{CL}, that is,
$$u(x)=\ln \Big[\frac{32 \lambda^2}{(4+\lambda^2 |x-x_0|^2)^2} \Big] \;\;\; \mbox{with $\lambda>0$, $x_0 \in \R^2$}.$$
Moreover, Dancer \& Farina \cite{DF} proved that \eqref{new1.1} with $\alpha=0$ admits classical stable at infinity solutions
if and only if $N=2$ or $N \geq 10$. Later on, Wang \& Ye \cite{WY} considered the nonautonomous case $\alpha \neq 0$.  They obtained
the following results:
\begin{itemize}
\item[(a)] \eqref{new1.1} admits no weak solution for any domain $\Omega \subset \R^N$ containing 0 provided $\alpha \leq -2$.
\item[(b)] \eqref{new1.1} does not admit any stable solution provided that $\alpha>-2$ and $2 \leq N<10+4 \alpha$.
\item[(c)] \eqref{new1.1} does not admit stable at infinity solution provided that $\alpha>-2$ and $2<N<10+4 \alpha^-$, where
$\alpha^-=\min (\alpha,0)$.
\end{itemize}

Clearly, (b) is sharp, since for $N \geq 10+4 \alpha$ and $\alpha>-2$, \eqref{new1.1} possesses radial stable solutions in $\R^N$, especially a singular stable solution
\begin{equation*}
U(x)=-(2+\alpha) \ln |x|+\ln [(2+\alpha)(N-2)].
\end{equation*}
For the same reason, (c) is sharp for $-2<\alpha \leq 0$. However, the following problem remained open: {\it For $\alpha>0$, does \eqref{new1.1} admit any stable at infinity solution provided that $10 \leq N< 10+4 \alpha$?}

\medskip
We will present here a negative (hence optimal) answer to this question.

\begin{thm}
\label{t1.1-}
For $\alpha>0$, \eqref{new1.1} does not admit any stable at infinity solution provided that $10 \leq N< 10+4 \alpha$.
\end{thm}

Let's look now the biharmonic equation \eqref{1.1-}. Let $u\in C(\R^N)$ be a weak solution to $\Delta^2 u=|x|^{\alpha} e^u$ in ${\mathcal D}'(\R^N)$. $u$ is said to be stable on $\O \subset \R^N$ if
\begin{align*}
\label{stable}
\int_{\O} |\Delta \phi|^2 dx - \int_{\O} |x|^{\alpha} e^u \phi^2 dx \geq 0,\quad \forall ~\phi\in C^2_c (\O).
\end{align*}
Similarly, a solution $u$ of \eqref{1.1-} is said stable at infinity, if $u$ is stable on $\R^N\backslash \mathcal {K}$ for some compact set $\mathcal {K}$; $u$ is said stable if it's stable in whole $\R^N$. We can define also the associated Morse index, again any finite Morse index solution $u$ on $\R^N$ is stable at infinity.

\medskip
When $\alpha = 0$, the properties of entire solutions of the equations as \eqref{1.1-} were studied by many authors, especially the study of stability of radial solutions of \eqref{1.1-} is well done, see \cite{AGG, BFFG, DGGW, W}. We will consider the radial solution for the non-autonomous case $\alpha>-2$ through the following initial value problem:
\begin{equation}
\label{01}
\begin{cases}
\begin{aligned}
&\Delta^2 u=|x|^{\alpha} e^u~~~~~~~~~\text{for}~~r\in [0,~R_{\delta,\beta})\\
&u'(0)=u'''(0)=0,\\
&u(0)=\delta,~~~\Delta u(0)=\beta,
\end{aligned}
\end{cases}
\end{equation}
where $R_{\delta,\beta}$ is the maximal interval of existence for given $\delta,~\beta\in\mathbb{R}$.

\medskip
Note that the equation \eqref{1.1-} is invariant under the scaling transformation:
$$
u_{\lambda}(x)=u(\lambda x)+(4+\alpha) \ln \lambda, ~\lambda>0.
$$
Remark also that the above transformation does not change the stability at infinity nor the stability of the solution. So we need only to consider $\delta =0$.

\medskip
For the existence of entire radial solutions, we have
\begin{prop}
\label{l1.4}
Let $\alpha>-2$ and $\delta = 0$. For $N\leq 2$, the initial value problem \eqref{01} admits no entire solution; for $N\geq 3$, there exists $\beta_0 < 0$ such that $R_\beta = \infty$ if and only if $\beta \leq \beta_0$.
\end{prop}

The following  result gives the asymptotic behavior of the border line entire solution $u_{\beta_0}$, called also {\it separatrix}.
\begin{thm}
\label{thm1.1}
Let $\alpha>-2$, $\delta = 0$ and $\beta_0$ be in Proposition \ref{l1.4}. Then the solution $u_{\beta_0}$ verifies:

(i) For $N=3$, as $r\to \infty$, $u_{\beta_0}(r) =a_1 r +a_2 +a_3 r^{-1}+O(e^{-cr})$ with $c>0$ and
$$a_1=-\frac{1}{8\pi}\int_{\mathbb{R}^3} |x|^{\alpha} e^{u_{\beta_0}} dx,~~~ a_2=\frac{1}{8\pi}\int_{\mathbb{R}^3} |x|^{1+\alpha} e^{u_{\beta_0}} dx,~~~ a_3=-\frac{1}{24\pi}\int_{\mathbb{R}^3} |x|^{2+\alpha} e^{u_{\beta_0}} dx.
$$

(ii) For $N=4$,
$$u_{\beta_0}(r)+(8+2\alpha)\ln r \to \frac{1}{8\pi^2}\int_{\mathbb{R}^4}  e^{u_{\beta_0}}|x|^{\alpha}\ln|x| dx \;\;\; \mbox{as $r\to +\infty$}.$$

(iii) For $N\geq 5$,
$$
u_{\beta_0}(r)+(4+\alpha) \ln r \to \ln\big[2(4+\alpha)(N-2)(N-4)\big] \;\;\; \mbox{as $r \to +\infty$}.
$$
\end{thm}

To state our results on the stability properties of entire radial solutions, we denote $N_{\alpha}$ the unique root in $(5, \infty)$ to the equation $$\frac{N^2(N-4)^2}{16}=2(4+\alpha)(N-2)(N-4).$$
Indeed, let $f(s) = s^2(s-4) - 32(4+\alpha)(s-2).$ It follows from $\alpha > -2$ that $f(5) < 0$ and a unique root exists in $(5, \infty)$ since $\lim_{s\to\infty} f(s)= \infty$ and $f$ is convex in $[5, \infty)$.

\begin{thm}
\label{thm1.2}
Let $\alpha>-2$, $\delta = 0$ and $u_\beta$  be a radial entire solution to \eqref{01} (hence $\beta \leq \beta_0$).

(i) For $N=3$ and 4, $u_{\beta}$ is unstable but stable at infinity for any $\beta\leq \beta_0$.

(ii) For $5\leq N<N_{\alpha}$, $u_{\beta}$ is stable at infinity for every $\beta<\beta_0$, while $u_{\beta_0}$ is unstable outside any compact set. Moreover, there exists $\beta_1<\beta_0<0$ such that $u_{\beta}$ is stable if and only if $\beta \leq \beta_1$.

(iii) If $N\geq N_{\alpha}$, $u_{\beta}$ is stable for every $\beta\leq \beta_0$.
\end{thm}

Therefore, we get a complete picture of the stability and/or stability at infinity of entire radial solutions to \eqref{1.1-} with $\alpha > -2$. A natural question is to understand general polyharmonic equations $(-\Delta)^m u = |x|^\alpha e^u$. However, the cases $m \geq 3$ present many new phenomena, where the situation is not completely understood, even for radial solutions with $\alpha = 0$, see for instance \cite{FF, HY}.

\medskip
Throughout this paper, $B_R$ denotes the open ball of radius $R$ centered at $0$. The constant $C>0$ denotes a generic number which may be different from line to line. Theorem \ref{t1.1-} is proved in section 2. In section 3, we show Proposition \ref{l1.4} and some basic estimates of entire solutions of \eqref{1.1-}. In sections 4-5, we present the proof of Theorems \ref{thm1.1} and \ref{thm1.2} respectively.

\section{Stability at infinity for solutions of \eqref{new1.1}}
\reset
We begin with an integral estimate which is a little variant of Proposition 2.2 of \cite{WY}, which is obtained by the stability at infinity with test function of the form $e^{\gamma u}\psi^m$.
\begin{prop}
\label{p1.1}
Let $\Omega$ be a domain (bounded or not) in $\R^N$. Let $u$ be a stable solution of \eqref{new1.1} on $\Om$ with $\alpha>-2$. Then
for any integer $m \geq 5$ and any $\gamma \in (0,2)$, there exists $C>0$ depending only on $m, \alpha$ and $\gamma$ such that
\begin{equation}
\label{1.4-1}
\int_\Omega \Big[|\nabla (e^{\gamma u})|^2+|x|^\alpha e^{(2 \gamma+1) u} \Big] \psi^{2m} \leq C \int_\Omega |x|^{-2 \gamma \alpha} (|\nabla \psi|^2+|\psi| |\Delta \psi|)^{2\gamma+1} dx,
\end{equation}
for all functions $\psi \in C_0^\infty (\Omega)$ verifying $\|\psi\|_\infty \leq 1$ in $\Omega$.
\end{prop}

Let $3 \leq N < 10 + 4\alpha$ and $u$ be a solution of \eqref{new1.1} stable outside a closed ball $\overline{B}_{R_*}$ with $R_* > 0$. Applying Proposition \ref{p1.1}, as for the estimate (2.8) in \cite{WY}, there exists a constant $C > 0$ such that for any $\gamma \in (0,2)$ and any $R> 8R_*$, it holds
\begin{align}
\label{new0}
\int_{\frac{R}{2} < |x|< R} \Big[|\nabla(e^{\gamma u})|^2 + |x|^{\alpha} e^{(2\gamma+1)u}\Big] dx \leq C R^{N-2(2\gamma+1)-2\gamma\alpha}.
\end{align}

Recall that here we study general weak solutions to \eqref{new1.1} without radial symmetry assumption. Let $\overline f$ denote the spherical average of function, i.e.
\begin{align}
\label{average}
\overline f(r) = \frac{1}{\omega_N}\int_{\mathbb{S}^{N-1}} f(r,\theta) d\theta, \quad \mbox{where } \omega_N = |\mathbb{S}^{N-1}|.
\end{align}
We have
\begin{itemize}
\item[$(i)$] Either $\lim_{r\to\infty} r^{2+\alpha}\overline{e^{u(r)}} =0$;
\item[$(ii)$] Or  there exist a constant $C_0 > 0$ and a sequence $r_k$ tending to $\infty$ such that
\begin{align}
\label{new-2}
r_k^{2+\alpha}\overline{e^{u(r_k)}} \geq C_0 > 0, \quad \forall\; k \in \N.
\end{align}
\end{itemize}

Suppose that $(i)$ holds true. Let $w=\overline u$ and $M > 0$ satisfying $\alpha+2-\frac{2}{(N-2)M} > 0$. For $r$ large enough, $-\Delta w(r) = r^\alpha \overline{e^u} \leq M^{-1}r^{-2}$. Integrating twice this inequality, similar to the end of the proof for Theorem 1.5 in \cite{WY}, we can claim that
\begin{align*}
r^{2+\alpha}\overline{e^u}(r) \geq r^{2+\alpha} e^{w(r)} \geq Cr^{\alpha+2-\frac{2}{(N-2)M}} \to \infty,
\end{align*}
hence contradicts $(i)$, hence the case $(i)$ cannot occur. Assume now we are in the case $(ii)$. Consider the scaling sequence
\begin{align*}
v_k(x) = u(r_kx) + (2 + \alpha)\ln r_k.
\end{align*}
Clearly, $v_k$ are weak solutions of \eqref{new1.1}. For $R_0 > 0$, since
\begin{align*}
& \int_{\frac{R_0}{2} < |x|< R_0} \Big[|\nabla(e^{\gamma v_k})|^2 + |x|^{\alpha} e^{(2\gamma+1)v_k}\Big] dx\\
 = &\; r_k^{2\gamma(2+\alpha)+2-N}\int_{\frac{r_kR_0}{2} < |y|< r_kR_0} \Big[|\nabla(e^{\gamma u})|^2 + |y|^{\alpha} e^{(2\gamma+1)u}\Big] dy,
\end{align*}
from \eqref{new0} we have
\begin{lem}
\label{lnew1}
Given any $R_0 > 0$ and $\gamma \in (0, 2)$, for $k$ large enough, there holds
\begin{align*}
\int_{\frac{R_0}{2} < |x|< R_0} \Big[|\nabla(e^{\gamma v_k})|^2 + |x|^{\alpha} e^{(2\gamma+1)v_k}\Big] dx \leq C < \infty.
\end{align*}
\end{lem}
Consequently, for any $R_0 > 0$, $(e^{v_k})$ is bounded in $H^1(A_0)$ where $A_0 := B_{R_0}\backslash \overline{B}_{R_0/2}$. Let $w_k(s) = \overline{e^{v_k}}(s)$, as $w_k'(s) = \overline{e^{v_k}\p_rv_k}(s)$, we check easily that as function of $s$, $\{w_k\}$ is bounded in $H^1(I_0)$ where $I_0 = [R_0/2, R_0]$. The Sobolev embedding implies that, up to a subsequence, $w_k$ convergence to $w$ in $C(I_0)$. Choosing now $R_0 = \frac{3}{2}$, by \eqref{new-2}, there holds
\begin{align}
\label{new4}
w_k(1) \geq C_0 > 0, \quad \mbox{hence } w(1) \geq C_0.
\end{align}
However, we have another estimate which leads to a contradiction.
\begin{lem}
\label{lnew2}
Let $3 \leq N < 10 + 4\alpha$, $u$ be a solution of \eqref{new1.1}, stable at infinity. Then
\begin{align*}
\lim_{R\to\infty}R^{2-N}\int_{B_R} |x|^\alpha e^u dx = 0.
\end{align*}
\end{lem}

Admitting the above result, noting that for any fixed $R_0 > 0$ and $r_k \to \infty$, we have
\begin{align*}
\omega_N \int_{I_0}s^{N-1+\alpha}w_k(s) ds = \int_{A_0} |y|^\alpha e^{v_k} dy & = r_k^{2-N}\int_{\frac{r_kR_0}{2} < |x|< r_kR_0} |x|^\alpha e^u dx \\
& \leq r_k^{2-N}\int_{B_{r_kR_0}} |x|^\alpha e^u dx \to 0, \quad \mbox{as } k\to \infty.
\end{align*}
Passing to the limit, as $w \geq 0$, we see that $w \equiv 0$ in $I_0$ which contradicts \eqref{new4}.

\medskip
In both case $(i)$ and $(ii)$, we prove the nonexistence of stable at infinity solution $u$ to \eqref{new1.1} when $3 \leq N < 10 + 4\alpha$. It remains to check Lemma \ref{lnew2}. Applying Proposition \ref{p1.1} with suitable cut-off function $\psi$, as in the proof of estimate (2.7) in \cite{WY}, we get that for any $\gamma \in (0, 2)$, there exist $C_1, C_2 > 0$ satisfying
\begin{align}
\label{wy1}
  \int_{2R_* < |x|< R} \Big[|\nabla(e^{\gamma u})|^2 + |x|^{\alpha} e^{(2\gamma+1)u}\Big] dx \leq C_1 + C_2 R^{N-2(2\gamma+1)-2\gamma\alpha}, \;\;\forall\; R> 8R_*.
  \end{align}
Let $\gamma_0$ satisfy $N-2(2\gamma_0+1)-2\gamma_0\alpha = 0$. We can check that $\gamma_0 \in (0, 2)$ since $3 \leq N < 10 + 4\alpha$. Therefore, tending $R$ to $\infty$ in \eqref{wy1} with $\gamma = \gamma_0$, we get
\begin{align}
\label{new5}
\int_{\R^N\backslash B_{2R_*}}|x|^{\alpha} e^{(2\gamma_0 +1)u} dx < \infty.
\end{align}
Moreover, by H\"older inequality, for $R > R_1 > 2R_*$,
\begin{align*}
\int_{B_R\backslash B_{R_1}} |x|^\alpha e^u dx & \leq \left(\int_{B_R\backslash B_{R_1}} |x|^\alpha e^{(2\gamma_0 +1)u} dx\right)^\frac{1}{2\gamma_0+1}\left(\int_{B_R\backslash B_{R_1}} |x|^\alpha dx\right)^\frac{2\gamma_0}{2\gamma_0+1}\\
& \leq \left(\int_{B_R\backslash B_{R_1}} |x|^\alpha e^{(2\gamma_0 +1)u} dx\right)^\frac{1}{2\gamma_0+1}\times C_{N, \alpha, \gamma_0}R^\frac{2\gamma_0(N+\alpha)}{2\gamma_0+1}\\
& = C_{N, \alpha, \gamma_0}\left(\int_{B_R\backslash B_{R_1}} |x|^\alpha e^{(2\gamma_0 +1)u} dx\right)^\frac{1}{2\gamma_0+1}R^{N-2}.
\end{align*}
Here $C_{N, \alpha, \gamma_0}$ is a constant independent on $R$ and $R_1$. We obtain then
\begin{align*}
R^{2-N}\int_{B_R} |x|^\alpha e^u dx \leq R^{2-N}\int_{B_{R_1}} |x|^\alpha e^u dx + C_{N, \alpha, \gamma_0}\left(\int_{B_R\backslash B_{R_1}} |x|^\alpha e^{(2\gamma_0 +1)u} dx\right)^\frac{1}{2\gamma_0+1}.
\end{align*}
Sending $R$ to $\infty$, there holds
\begin{align*}
\limsup_{R\to \infty} \left(R^{2-N}\int_{B_R} |x|^\alpha e^u dx\right) \leq C_{N, \alpha, \gamma_0}\left(\int_{\R^N\backslash B_{R_1}} |x|^\alpha e^{(2\gamma_0 +1)u} dx\right)^\frac{1}{2\gamma_0+1}.
\end{align*}
Tending now $R_1$ to $\infty$, we get the claim by \eqref{new5}. So we are done.\qed

\section{Existence of entire radial solutions to \eqref{1.1-}}
\reset
Here we study the existence of global radial solutions to \eqref{01}, where we will use the following well known comparison result.
\begin{lem}
\label{l2.1}
Let $f=f(r, s): \R \times \R \to \R$ be nondecreasing in $s$. Let $u, v\in C^{2m}([\rho, R))$ be two radial functions such that $\Delta^m u - f(r, u)\geq \Delta^m v -f(r, v)$ in $[\rho, R)$ for some $m\geq 1$ and $R > \rho\geq 0$. Suppose that
\begin{align}
\Delta^ku(\rho) \geq \Delta^kv(\rho), \;\; (\Delta^ku)'(\rho) \geq (\Delta^kv)'(\rho) \quad \forall\; 0\leq k \leq m-1.
\end{align}
Then for any $r\in[\rho, R)$ and any $0 \leq k \leq m-1$, there hold
$$
\Delta^ku(r) \geq \Delta^kv(r),\;\;( \Delta^ku)'(r) \geq (\Delta^kv)'(r).
$$
\end{lem}

\noindent
{\bf Proof of Proposition \ref{l1.4}.}
Fix $\delta = 0$. Let $u$ be a continuous radial function, define
$$w(r) := T_{\beta, \alpha}(u) := \beta + \frac{1}{N-2}\int^r_0 s^{1+\alpha} e^u ds -\frac{r^{2-N}}{N-2}\int^r_0 s^{N-1+\alpha} e^u ds.$$
Obviously $\Delta w(r)= r^{\alpha}e^u$ and
$$\|w-\beta\|_{L^{\infty}([0,r])}\leq \frac{1}{N-2}\int^r_0 s^{1+\alpha} e^u ds,\quad \forall \; r > 0.$$
To get the local existence of solution to \eqref{01}, we can apply just the usual contraction argument with $\phi(u) = T_{0, 0}\circ T_{\beta, \alpha}(u)$.

\medskip
Suppose that $R_\beta = \infty$. Since
 \begin{equation}
 \label{new-1.3}
 \begin{aligned}
 \Delta u(r)&=\Delta u(0)+\int^r_0 s^{1-N} \int^s_0 \sigma^{N-1} \Delta^2 u(\sigma) d\sigma ds\\
 &= \beta +\int^r_0 s^{1-N} \int^s_0 \sigma^{N+\alpha-1}  e^{u(\sigma)} d\sigma ds,
 \end{aligned}
\end{equation}
$\Delta u$ is increasing in $r$. We claim that $\lim_{r\to\infty} \Delta u_\beta (r) = \ell \leq 0$ and $\beta \leq 0$. Otherwise, $\Delta u\geq C > 0$ at infinity, then $u \geq Cr^2$ and $r^\alpha e^{u(r)} \geq Cu^2(r)$ for large $r$. Seeing the following Mitidieri \& Pohozaev's result (see \cite{MP}), we reach a contradiction with $R_\beta = \infty$.
\begin{lem}
\label{2.4-}
Let $f\in C(\R\times\R)$. Assume that there exists $p > 1$ such that $f(r,s) \geq s^p$ for large $s$. Let $u$ be a solution to $\sum _{1\leq k \leq m} a_k u^{(k)}(t) = f(t,u(t))$ in $[t_0, \infty)$, with $a_k\in \R$ and $m\geq 1$. Then we cannot have $\lim_{t\to\infty} u(t)= +\infty$.
\end{lem}

We conclude then: If $u_\beta$ is a entire solution to \eqref{01}, there hold $\beta < 0$ and $\lim_{r\to\infty} \Delta u_\beta(r) \leq 0$. Consequently $(\D u_\beta)(r) < 0$ for $r>0$ and $u_\beta$ is decreasing in $(0, \infty)$.

\medskip
An immediate consequence is the nonexistence result in lower dimensions. As $r^{N-1} (\Delta u)'(r)$ is increasing in $r$ and $r^{N-1} (\Delta u)'(r)>0$ for $r \in (0, \infty)$, if $u$ is an entire solution to \eqref{01} with $N\leq 2$, we get $\lim_{r\to\infty}\D u(r) =\infty$, which is impossible, hence \eqref{1.1-} has no entire solutions for $N \leq 2$.

\medskip
Fix now $N\geq 3$ and $\beta<0$. The comparison principle with the operator $L (\phi):=\Delta^2 \phi -|x|^\alpha e^{\phi}$ for radial fucntions (see Lemma \ref{l2.1}) implies
\begin{equation}
\label{1.3}
u_\beta (r) \geq \frac{\beta}{2N} r^2, \;\;\; \mbox{for $r \in (0, R_\beta]$}.
\end{equation}
Moreover, let $u_1(r) =\gamma_1 r^{4+\alpha}-Br^2-A\ln(1+r^2)$, then
$$u_1'(r)=(4+\alpha) \gamma_1 r^{3+\alpha}-2Br -\frac{2Ar}{1+r^2},$$
$$\Delta {u_1}(r)=(4+\alpha)(N+2+\alpha)\gamma_1 r^{2+\alpha}-2NB -\frac{2A(N-2)}{1+r^2}-\frac{4A}{(1+r^2)^2},$$
and
$$
(\Delta{u_1})'(r)=(4+\alpha)(2+\alpha)(N+2+\alpha) \gamma_1 r^{1+\alpha}+\frac{4A(N-2)r}{(1+r^2)^2}+\frac{16Ar}{(1+r^2)^3}.
$$
Fix $\gamma_1 > 0$ such that $\D^2(\gamma_1 r^{4+\alpha}) = 2r^\alpha$, we have
\begin{align*}
\Delta ^2u_1(r)= r^{\alpha}\left[ 2+\frac{4A(N-2)(N-4)r^{-\alpha}}{(1+r^2)^2}+\frac{32A(N-4)r^{-\alpha}}{(1+r^2)^3}+\frac{96A r^{-\alpha}}{(1+r^2)^4}\right]=:r^{\alpha}H(r).
\end{align*}
For $N \geq 3$ and $r \in (0, 1]$, there holds
$$
H(r)\geq 2+ \frac{16A r^{-\alpha}}{(1+r^2)^4} \geq 2 + Ar^{-\alpha}.
$$
Take now $A >0$ and a suitably large $B>0$ (to be chosen later) such that
$$2 + Ar^{-\alpha} \geq e^{\gamma_1 r^{4+\alpha}} \geq \frac{e^{\gamma_1 r^{4+\alpha}}}{e^{B r^2}(1+r^2)^A} = e^{u_1(r)}, \quad \forall\; r \in (0, 1].$$
Hence $u_1$ is a supersolution to \eqref{01} in $B_1$ with $u_1(0) = 0$, $\Delta u_1(0) = -2N(A+B)$. By the comparison principle,
$$u(r)\leq u_1(r),\quad u'(r)\leq u_1'(r),\quad \Delta u(r)\leq \Delta u_1(r), \quad (\Delta u)'(r)\leq (\Delta u_1)'(r)\quad\forall\;0\leq r\leq 1,$$
provided $\beta \leq -2N(A+B)$. This and \eqref{1.3} imply $R_\beta > 1$. To prove $R_\beta = \infty$ for $\beta \to -\infty$, we divide our study into two cases.

\smallskip\noindent
{\it Case I: $N\geq 4$}. Define $u_2(r)=-\bar Br^2-\bar A\ln (1+r^2)$ with $\bar A>0$. Then $u_2$ satisfies
$$
\Delta^2{u_2}(r) =\frac{4\bar A(N-2)(N-4)}{(1+r^2)^2}+\frac{32\bar A(N-4)}{(1+r^2)^3}+\frac{96\bar A}{(1+r^2)^4}\geq\frac{96\bar A}{(1+r^2)^4}=:|x|^{\alpha} e^{u_2}\times \frac{1}{\Phi(r)},
$$
where
$$\Phi(r) = \frac{|x|^{\alpha}(1+r^2)^{4-\bar A}}{96\bar Ae^{\bar Br^2}}.$$
A simple calculation implies that there exists $\bar A > A$ such that $(\Delta u_1)'(1)=(\Delta u_2)'(1)$. Then if $\bar B = B/2$ with sufficiently large $B$, we get
\begin{equation}
\label{new1.3}
u_1(1)\leq u_2(1), \quad u_1'(1)\leq u_2'(1), \quad\Delta u_1(1)\leq \Delta u_2(1), \quad(\Delta u_1)'(1)\leq (\Delta u_2)'(1),
\end{equation}
and $\sup_{r\geq 1}\Phi\leq 1$. Clearly \eqref{new1.3} holds also for $u$. Using again the comparison principle, we get $u(r)\leq u_2(r)$ for $r \in [1, R_\beta)$. Hence no blow-up occurs, that is $R_\beta = \infty$ for $\beta \leq -2N(A+B)$ with sufficiently large $B$.

\smallskip\noindent
{\it Case II: $N=3$}. Define here $u_3(r)=-\bar Br^2+\bar{A}\ln (1+r)$ with $\bar{A}>0$. Then $u_3$ verifies
$$
u_3'(r)=-2\bar Br +\frac{\bar A}{1+r}, \quad u_3''(r)=-2\bar B-\frac{\bar{A}}{(1+r)^2},\quad u_3'''(r)=\frac{2\bar{A}}{(1+r)^3}
$$
and
$$
\Delta^2u_3(r)=\frac{2\bar A r^3}{(1+r)^3}+\frac{6\bar A r^3}{(1+r)^4}\geq \frac{2\bar A r^3}{(1+r)^3}.
$$
Taking $\bar A$ such that $u'''_1(1)=u'''_3(1)$ and $\bar B = B/2$, we obtain that, for sufficiently large $B>0$, $\Delta^2 u_3\geq r^{\alpha} e^{u_3}$ in $\R^3\backslash B_1$ and \eqref{new1.3} if we replace $u_2$ by $u_3$. We have the same conclusion as in {\it Case I}.

\medskip
Therefore,
\begin{align}
\label{b0}
\beta_0 = \sup\big\{\beta \in \R \;\mbox{such that}\; R_\beta = \infty\}
\end{align}
is well defined, $\beta_0 \leq 0$ and for any $\beta <\beta_0$, the solution $u_{\beta}$ is global. Furthermore, for $\beta<\beta_0$, by \eqref{new-1.3} and the comparison principle with $\Phi(r):=-Cr^2$, we get readily the low and upper bound for $u_\beta$
\begin{equation}
\label{1.3-}
u_{\beta} \geq \frac{\beta}{2N} r^2\quad \text{and}\quad u_{\beta}\leq -\frac{\beta_0-\beta}{2N} r^2\quad \text{for all}~r\in[0,\infty).
\end{equation}
So $\lim_{r\to +\infty} \Delta u_{\beta}(r) <0$, if $\beta<\beta_0$.

\medskip
It remains to prove that $R_{\beta_0} = \infty$ and $\beta_0 < 0$. Consider $u_{\beta}$ with $\beta<\beta_0$. By the monotonicity of $u_{\beta}$, there holds
$$
(\Delta u_{\beta})'(r)=\frac{1}{r^{N-1}} \int^r_0 s^{N-1+\alpha} e^{u_{\beta}(s)} ds\geq \frac{r^{1+\alpha}}{N+\alpha} e^{u_{\beta}(r)}.
$$
We deduce that
$$
\Delta u_{\beta} (r)\geq \beta + \frac{1}{N+\alpha} \int^r_0 s^{1+\alpha} e^{u_{\beta}(s)} ds\geq \beta+\frac{1}{(2+\alpha)(N+\alpha)}r^{2+\alpha} e^{u_{\beta}(r)}.
$$
Since $\Delta u_{\beta}<0$ in $\R^N$,
$$
e^{u_{\beta}(r)}\leq -\frac{\beta(2+\alpha)(N+\alpha)}{r^{2+\alpha}}\quad\forall\; r > 0.
$$
By the continuity of $u_{\beta}$ w.r.t the initial value, we can conclude that
\begin{align}
\label{new1}
e^{u_{\beta_0}(r)}\leq -\frac{\beta_0(2+\alpha)(N+\alpha)}{r^{2+\alpha}}\quad\forall \; 0<~r<R_{\beta_0}.
\end{align}
This yields that $\beta_0<0$ and there is no blow-up at finite value of $r$ for $u_{\beta_0}$, so $R_{\beta_0}=+\infty$. \qed

\medskip
We now give a refined global estimate (comparing to \eqref{new1}) for the entire solution with $\beta\leq \beta_0$.
\begin{lem}
\label{l2.2}
Let $N\geq 3$. Then there exists $C_{N, \alpha}>0$ such that for any $\beta \leq \beta_0$, $u_{\beta}$ verifies
\begin{equation}
\label{2.1}
u_\beta(r) \leq-(4+\alpha) \ln r + C_{N, \alpha},\;\;\forall ~r>0.
\end{equation}
\end{lem}

\noindent
{\bf Proof.} For any $r_0 > 0$, let $M_0 :=u_{\beta}(r_0)+(4+\alpha)\ln r_0$. Define $v(x):= u_{\beta}(r_0|x|)+(4+\alpha)\ln r_0 -M_0,$
then $v$ satisfies the equation
\begin{equation}
\label{2.2}
\left \{ \begin{array}{ll} \Delta^2 v= \lambda |x|^{\alpha} e^{v} \geq 0 \;\; & \mbox{in $B_1$},\\
v=0, \;\; -\Delta v= -r_0^2 \Delta u_\beta(r_0) > 0\;\; & \mbox{on $\partial B_1$},
\end{array} \right.
\end{equation}
where $\lambda =e^{M_0}$. This implies that $v > 0$ in $B_1$ and $v$ is a supersolution to the following Navier boundary value problem
\begin{equation}
\label{2.3}
\left \{\begin{array}{ll} \Delta^2 w = \lambda |x|^{\alpha} e^{w} \;\; & \mbox{in $B_1$},\\
w= \Delta w =0 \;\; & \mbox{on $\partial B_1$}.
\end{array} \right.
\end{equation}
Hence \eqref{2.3} admits a positive solution by the monotone iteration method. On the other hand, since $|x|^{\alpha}\in L^q(B_1)$ for some $q > \frac{N}{4}$, letting $\lambda^*$ be the first eigenvalue given by
$$
\lambda^* = \inf\left\{\int_{B_1} |\Delta \phi|^2 dx, \phi\in H^2(B_1)\cap H^1_0(B_1), \int_{B_1} |x|^{\alpha}\phi^2dx = 1\right\},
$$
it's well-known that $\l^* \in (0, \infty)$, and there exists a minimizer $\phi$ satisfying $\phi > 0$ in $B_1$ and
\begin{equation}
\label{2.4}
\left \{ \begin{array}{ll} \Delta^2 \phi = \lambda^* |x|^{\alpha}\phi \;\; & \mbox{in $B_1$},\\
\phi= \D \phi =0 \;\; & \mbox{on $\partial B_1$}.
\end{array} \right.
\end{equation}
Multiplying $\phi$ and $w$ on both sides of \eqref{2.3} and \eqref{2.4} respectively and  integrating by parts on $B_1$, we deduce that
\begin{align*}
0 = \int_{B_1} |x|^\alpha\big(\l e^{w} - \l^* w\big)\phi dx.
\end{align*}
This is possible only if $\lambda <\lambda^*$, since otherwise $\l e^{w} - \l^* w > 0$ in $B_1$. Hence $u_\beta(r_0) + (4+\alpha)\ln r_0 < \ln\l^*$. Since $r_0 > 0$ is arbitrary, \eqref{2.1} holds true with $C_{N,\alpha}= \ln \l^*$. \qed

\begin{prop}
\label{newp1}
For any $\beta \leq \beta_0$, $\lim_{r\to\infty} \Delta u_{\beta} \leq 0$ and $\lim_{r\to\infty} \Delta u_{\beta} = 0$ if and only if $\beta =\beta_0$.
\end{prop}

\noindent
{\bf Proof.} Notice that if $\beta<\beta_0$, $\lim_{r\to\infty} \D u_{\beta}(r)< 0$. We only need to show that
$$\mbox{if $\lim_{r\to\infty} \D u_{\beta}(r) = \ell < 0$, then $\beta < \beta_0$.}$$
For simplicity, we omit the index $\beta$. We can claim then
$$u(r)\sim \frac{\ell}{2N} r^2,\quad u'(r)\sim \frac{\ell}{N}r \quad \text{as}~ r\to \infty$$ and $$(\Delta u)'(r)=r^{1-N}\int^r_0 s^{N-1+\alpha} e^u ds.$$
For radial functions $f$ and $g$, by ``$f\sim g$ near $\infty$" we mean $\lim_{r\to\infty}\frac{f(r)}{g(r)}=1$.

\medskip
For $N\geq 4$, consider $\tilde{u}(r) =-\epsilon r^2 -A \ln (1+r^2)$, where $\epsilon$ and $A$ are positive constants to be determined. Direct computations show that
$$
\tilde{u}'(r)=-2\epsilon r -\frac{2Ar}{1+r^2}, \Delta \tilde{u}(r) =-2N\epsilon -\frac{2A(N-2)}{1+r^2}-\frac{4A}{(1+r^2)^2},
$$
$$(\Delta \tilde{u})'(r)=\frac{4A(N-2)r}{(1+r^2)^2}+\frac{16Ar}{(1+r^2)^3}$$
and
$$
\Delta^2\tilde{u}(r) =\frac{4A(N-2)(N-4)}{(1+r^2)^2}+\frac{32A(N-4)}{(1+r^2)^3}+\frac{96A}{(1+r^2)^4}.
$$
Fix
$$0 < \epsilon < -\frac{\ell}{2N}, \quad A> \frac{1}{2N}\int^{\infty}_0 r^{N-1+\alpha} e^u dx.$$
There exists $r_0$ large such that for any $r \geq r_0$, $\Delta^2 \tilde u>r^{\alpha} e^{\tilde{u}}$,
\begin{equation}
\label{2.5}
(\Delta u)'(r)<(\Delta \tilde{u})'(r),\quad \Delta u(r)<\Delta\tilde{u}(r),\quad u'(r)<\tilde{u}'(r)\quad\mbox {and} \quad u(r)<\tilde{u}(r).
\end{equation}
The continuity with respect to the initial data implies that for $\beta' > \beta$ but close enough, the estimates \eqref{2.5} holds true for $u_{\beta'}$ and $r= r_0$. Applying once more the comparison principle, we see that $u_{\beta'} < \tilde u$ for all $r_0 < r < R_{\beta'}$. Therefore $u_{\beta'}$ will never blow up, so it is global. By the definition of $\beta_0$, we get $\beta' \leq \beta_0$, so $\beta < \beta_0$.

\medskip
For $N=3$, the equation $\Delta^2 u=r^{\alpha} e^u$ reads $(r^4u''')'=r^{4+\alpha}e^u$. Consider the function $\widehat{u}$ defined by $\widehat u(r) :=-\epsilon r^2+A\ln(1+r).$
As the above, with $\epsilon > 0$ small enough and $A$ large enough, we can check that there exists $r_0 > 1$ such that $\Delta^2 \widehat u \geq r^{\alpha} e^{\widehat{u}}$ and \eqref{2.5} hold true for $r \geq r_0$ with $\widehat u$ instead of $\tilde u$. We can conclude very similarly as above.\qed

\section{Asymptotic behaviors of Separatrix radial solutions to \eqref{1.1-}}
\reset
Here we will prove Theorem \ref{thm1.1}. For simplicity, we omit the index $\beta_0$ and $u$ stands always $u_{\beta_0}$ in this section. We consider separately the case $N = 3$, $N = 4$ and $N \geq 5$

\subsection{\it Case 1: $N=3$} Thanks to \eqref{2.1}, $r^{2+\alpha}e^{u(r)} \in L^1(\R_+)$, there holds then
$$
\Delta u(r) = -\int^{\infty}_r s^{-2}\int_0^s \sigma^{2+\alpha} e^{u(\sigma)} d\sigma ds \sim -\frac{1}{r}\int_0^\infty \sigma^{2+\alpha}e^{u(\sigma)} d\sigma,\quad \mbox{as } r\to\infty.
$$
Let
$$
a :=-\int_0^\infty \sigma^{2+\alpha}e^{u(\sigma)} d\sigma =-\frac{1}{4\pi }\int_{\mathbb{R}^3} |x|^{\alpha} e^{u(x)} dx<0.
$$
Then we have
$$u'(r) = r^{-2}\int^{r}_0\Delta u(\sigma) \sigma^2 d\sigma \to \frac{a}{2}, \quad \mbox{as } r\to\infty.$$
So $u(r) \leq -Cr$ for $r$ large. More precisely,
$$
\begin{aligned}
-\Delta u(r)&=\int^{\infty}_r s^{-2}\int^s_0 \sigma ^{2+\alpha} e^{u(\sigma)} d\sigma ds\\
&=\int^r_0\sigma ^{2+\alpha} e^{u(\sigma)} \int^{\infty}_r s^{-2} ds d\sigma + \int^{\infty}_r\sigma ^{2+\alpha} e^{u(\sigma)} \int^{\infty}_{\sigma} s^{-2} ds d\sigma \\
&=-\frac{a}{r} -\frac{1}{r}\int^{\infty}_r \sigma ^{2+\alpha}e^{u(\sigma)} d\sigma + \int^{\infty}_r \sigma^{1+\alpha} e^{u(\sigma)} d\sigma.
\end{aligned}
$$
Therefore
$$
(r^2u')' =ar+r\int^{\infty}_r \sigma ^{2+\alpha}e^{u(\sigma)} d\sigma -r^2\int^{\infty}_r \sigma^{1+\alpha} e^u d\sigma.
$$
Integrating over $[0,~r]$, there holds
$$
\begin{aligned}
r^2 u'(r)&=\frac{a r^2}{2}+\int^r_0s\int^{\infty}_s \sigma ^{2+\alpha}e^{u(\sigma)} d\sigma -\int^r_0 s^2\int^{\infty}_s \sigma^{1+\alpha} e^{u(\sigma)} d\sigma\\
&=\frac{a r^2}{2}+\int^\infty_0 \sigma ^{2+\alpha}e^{u(\sigma)} \int^{\min(\sigma, r)}_0 s ds d\sigma -\int^\infty_0\sigma^{1+\alpha} e^{u(\sigma)} \int^{\min(\sigma, r)}_0 s^2ds d\sigma\\
&=\frac{ar^2}{2} +\frac{1}{6}\int^r_0 \sigma^{4+\sigma} e^{u(\sigma)} d\sigma +\frac{r^2}{2} \int^{\infty}_r \sigma^{2+\alpha}e^{u(\sigma)} d\sigma -\frac{r^3}{3}\int^{\infty}_{r} \sigma^{1+\alpha} e^{u(\sigma)} d\sigma.
\end{aligned}
$$
We deduce then
\begin{align*}
u(r) & =\frac{ar}{2} +\frac{1}{6}\int^r_0 \frac{1}{s^2}\int^s_0 \sigma^{4+\sigma} e^{u(\sigma)} d\sigma +\frac{1}{2} \int_0^r\int^{\infty}_s \sigma^{2+\alpha}e^{u(\sigma)}d\sigma -\int^r_0\frac{s}{3}\int^{\infty}_{s} \sigma^{1+\alpha} e^{u(\sigma)} d\sigma\\
& =\frac{ar}{2}+\frac{1}{2}\int^r_0 \sigma^{3+\alpha}e^{u(\sigma)} d\sigma - \frac{1}{6r}\int^r_0\sigma^{4+\alpha} e^{u(\sigma)} d\sigma - \frac{r^2}{6}\int^{\infty}_r \sigma^{1+\alpha}e^{u(\sigma)} d\sigma.
\end{align*}
Now it is easy to get the claimed expansion for $u(r)$ at infinity.

\subsection{\it Case 2: $N=4$} First, we claim $$u(r) =-(8+2\alpha)\ln r  + \frac{1}{4 \omega_3}\int_{\mathbb{R}^4} (\ln r) r^{\alpha} e^u dx +o(1)~~~\quad \text{as}~~~r\to +\infty.$$
Indeed, as
$$\int^s_0 \frac{\sigma^3}{1+\sigma^4} d\sigma=\frac{1}{4}\ln(s^4+1)\leq 1+\ln s, \quad\mbox{for $s\geq 1$},$$
using \eqref{2.1}, we get
$$
|\Delta u(r)| \leq C\int^{\infty}_r s^{-3} (1+\ln s ) ds \leq C(r^{-2}+r^{-2} \ln r),~~~~~~~\text {for}~~r\geq 1.
$$
Consequently
$$
\begin{aligned}
|u(r)|&\leq C+C \int^r_1 s^{-3} \int^s_1 (\sigma^{-2}+\sigma^{-2} \ln \sigma)\sigma^3 d\sigma ds\\
&\leq C+ C \int^r_1 s^{-3} (1+s^2 +s^2 \ln s)ds\\
&\leq C+C (\ln r +(\ln r)^2).
\end{aligned}
$$
Hence, $|u(r)|=O((\ln r)^2)$ as $r \to \infty$. More precisely, we claim that $u(r)=O(\ln r)$ as $r\to \infty$. Back to the equation $\Delta^2 u=|x|^{\alpha} e^u$ in $\mathbb{R}^4$. As $u$ is radial,
$$
u^{(4)} +\frac{6}{r} u^{(3)}+\frac{3}{r^2} u''-\frac{3}{r^3}u'=r^{\alpha} e^u.
$$
Let $w(s)=u(r)+(4+\alpha)s$ with $s = \ln r$, then $w$ satisfies
\begin{equation}\label{1.4}
w^{(4)}-4 w''=e^w,
\end{equation}
i.e. $(\frac{\partial}{\partial s}+2)(\frac{\partial}{\partial s}-2)w''=e^w>0$. Denote $z=(\frac{\partial}{\partial s}-2)w''$, we see $(e^{2s} z)'>0$.
Integrating over $(-\infty, s)$, as $u(0)=0=(\Delta u)'(0),~~\Delta u(0)=\beta_0$, so
$$
\begin{aligned}
\lim_{s\to -\infty} e^{2s} z&=\lim_{s\to -\infty} e^{2s}(w'''-2w'') =\lim_{r\to 0} r^2\big[u'''(r) r^3+u''(r) r^2-u'(r) r\big]=0.
\end{aligned}
$$
Hence $e^{2s}(\frac{\partial}{\partial s}-2)w''\geq 0,$ or equally $(e^{-2s} w'')'\geq 0.$
As $\lim_{r\to +\infty} \Delta u =0$, we deduce
\begin{equation}\label{06}
\lim_{s\to +\infty} e^{-2s}w''=\lim_{s\to +\infty} e^{-2s}(u''e^{2s}+ u'e^s)=\lim_{r\to +\infty} (u''+\frac{u'}{r})=0.
\end{equation}
Therefore $e^{-2s} w''\leq 0$,  i.e. $w$ is concave. Hence $\lim_{s\to +\infty} w= \ell \in \overline{\mathbb{R}}$ exists. Seeing Lemma \ref{2.4-}, we can claim $\ell =-\infty$.

\medskip
Let $\xi=w''$, then the equation \eqref{1.4} reads $\left(e^{2s}(\frac{\partial}{\partial s}-2) \xi\right)'=e^{2s} e^w$ and $\lim_{s\to-\infty} \xi(s)=\lim_{s\to-\infty} \xi'(s)=0$, hence
$$
(e^{-2s}\xi)'=e^{-4s}\int^s_{-\infty} e^{2t} e^w dt.
$$
By \eqref{06}, there holds
$$
\begin{aligned}
e^{-2s} \xi&=-\int^{+\infty}_s e^{-4t}\int^t_{-\infty} e^{2\sigma} e^w d\sigma dt=-\frac{e^{-4s}}{4} \int^s_{-\infty}e^{2\sigma} e^w d\sigma -\frac{1}{4} \int^{+\infty}_s e^{-2\sigma} e^w d\sigma.
\end{aligned}
$$
i.e.
$$
w''=-\frac{e^{-2s}}{4} \int^s_{-\infty}e^{2\sigma} e^w d\sigma -\frac{e^{2s}}{4} \int^{+\infty}_s e^{-2\sigma} e^w d\sigma.
$$
On the other hand, $\lim_{s\to -\infty} w'(s)=\lim_{r\to 0} u'(r) r+(4+\alpha)=4+\alpha$. So we obtain
$$
\begin{aligned}
&\;\; w'(s)-(4+\alpha)\\=&-\frac{1}{4}\int^s_{-\infty} e^{-2t} \int^t_{-\infty}e^{2\sigma} e^w d\sigma -\frac{1}{4} \int^s_{-\infty} e^{2t}\int^{+\infty}_t e^{-2\sigma} e^w d\sigma\\
=&\;-\frac{1}{4}\int^s_{-\infty}e^{2\sigma} e^w\int^s_{\sigma} e^{-2t} dt d\sigma - \frac{1}{4}\left[\int^s_{-\infty}e^{-2\sigma}e^w \int^{\sigma}_{-\infty} e^{2t} dt d\sigma +\int^{+\infty}_s e^{-2\sigma}e^w \int^s_{-\infty} e^{2t} dt d\sigma\right]\\
=&\;\frac{1}{8}e^{-2s}\int^s_{-\infty} e^{2\sigma} e^w d\sigma -\frac{1}{8}e^{2s} \int^{+\infty}_s e^{-2\sigma}e^w d\sigma -\frac{1}{4}\int^s_{-\infty} e^w d\sigma.
\end{aligned}
$$
Now integrating over $(0,~s)$, we have
$$
\begin{aligned}
&\;\; w(s)-w(0)\\=&(4+\alpha) s+ \frac{1}{8}\int^s_0 e^{-2t}\int^t_{-\infty} e^{2\sigma} e^w d\sigma dt -\frac{1}{8}\int^s_0 e^{2t} \int^{+\infty}_t e^{-2\sigma}e^w d\sigma dt -\frac{1}{4}\int^s_0 \int^t_{-\infty} e^w d\sigma dt\\
=&\;(4+\alpha)s+\frac{1}{8}\left[\int^0_{-\infty} e^{2\sigma} e^w \int^s_0 e^{-2t}dt d\sigma +\int^s_0 e^{2\sigma} e^w \int^s_{\sigma} e^{-2t} dt d\sigma\right]\\
&\;-\frac{1}{8}\left[ \int^s_0 e^{-2\sigma} e^w \int^{\sigma}_0 e^{2t} dt d\sigma +\int^{+\infty}_s e^{-2\sigma} e^w \int^s_0 e^{2t} dt d\sigma\right]-\frac{1}{4} s\int^s_{-\infty} e^w d\sigma+\frac{1}{4}\int^s_0 \sigma e^w d\sigma\\
=&\;(4+\alpha)s +\frac{1}{16}\int^0_{-\infty} e^{2\sigma} e^w d\sigma +\frac{1}{16}\int^{+\infty}_0 e^{-2\sigma}e^w d\sigma -\frac{1}{16} e^{-2s}\int^s_{-\infty} e^{2\sigma} e^w d\sigma \\ &\; -\frac{1}{16}\int^{+\infty}_s e^{2(s-\sigma)} e^w d\sigma -\frac{s}{4}\int^{+\infty}_{-\infty} e^w d\sigma +\frac{s}{4}\int^{+\infty}_s e^w d\sigma +\frac{1}{4}\int^s_0 \sigma e^w d\sigma.
\end{aligned}
$$
As $w''(s)\leq 0$ i.e. $w$ is concave and $\lim_{s\to +\infty} w(s)=-\infty$. So $\lim_{s\to +\infty} w'=\bar{\ell} \in [-\infty, 0)$, hence $ w(s)\leq -c s$ for some $c>0$ and $s$ large. This means $s^{\tau}e^ w\in L^1(\mathbb{R})$ for any $\tau \in \mathbb{R}$. Therefore $$u(r)=-\frac{1}{4 \omega_3}\int_{\mathbb{R}^4} |x|^{\alpha} e^u dx \ln r+ O(1) \quad \mbox{as $r\to +\infty.$}$$

Next we estimate the constant term (in short CT) in the above expansion of $u(r)$.
  $$
  \begin{aligned}
  CT&=w(0)+\frac{1}{16}\int^0_{-\infty} e^{2\sigma} e^w d\sigma +\frac{1}{16}\int^{+\infty}_0 e^{-2\sigma}e^w d\sigma  +\frac{1}{4} \int^{\infty}_0 \sigma e^{w(\sigma)} d\sigma\\
  &=u(1)+\frac{1}{16}\int^1_0 r^{5+\alpha} e^u dr +\frac{1}{16} \int^{\infty}_1 r^{1+\alpha} e^u dr + \frac{1}{4} \int^{\infty}_1  r^{3+\alpha} \ln r e^u dr.
  \end{aligned}
  $$
To estimate $u(1)$, recall the equation $\Delta^2 u= |x|^{\alpha} e^u$ in $\mathbb{R}^4$. As $\lim_{r\to +\infty} \Delta u =0$, there holds
$$
\begin{aligned}
-\Delta u(r) =\int^{\infty}_r s^{-3} \int^s_0\rho^{3+\alpha} e^u d\rho ds &=\int^r_0 \rho^{3+\alpha}e^u\int^{\infty}_r s^{-3}ds d\rho + \int^{\infty}_r \rho^{3+\alpha}e^u \int^{\infty}_{\rho} s^{-3} ds d\rho\\
&=\frac{r^{-2}}{2} \int^r_0 \rho^{3+\alpha} e^u d\rho +\frac{1}{2}\int^{\infty}_r \rho^{1+\alpha} e^u d\rho.
\end{aligned}
$$
Hence
$$
\begin{aligned}
-r^3 u'&=\frac{1}{2} \int^r_0 s \int^s_0\rho^{3+\alpha} e^u d\rho ds +\frac{1}{2} \int^r_0 s^3\int^{\infty}_s \rho^{1+\alpha} e^u d\rho ds\\
&=\frac{r^2}{4}\int^r_0 \rho^{3+\alpha} e^u d\rho - \frac{1}{8}\int^r_0 \rho^{5+\alpha} e^u d\rho +\frac{r^4}{8}\int^{\infty}_r \rho^{1+\alpha} e^u d\rho.
\end{aligned}
$$
Consequently, we get
$$
\begin{aligned}
u(1)=&-\frac{1}{4}\int^1_0\frac{1}{r}\int^r_0 \rho^{3+\alpha}d\rho dr +\frac{1}{8} \int^1_0\frac{1}{r^3} \int^r_0 \rho^{5+\alpha} e^u d\rho dr -\frac{1}{8} \int^1_0 r\int^{\infty}_r \rho^{1+\alpha}e^u d\rho dr\\
=&-\frac{1}{4}\int^1_0\rho^{3+\alpha} e^u\int^1_{\rho} r^{-1} dr d\rho + \frac{1}{8} \int^1_0 \rho ^{5+\alpha} e^u \int^1_{\rho} r^{-3} dr d\rho\\
&-\frac{1}{8} \int^1_0\rho^{1+\alpha}e^u \int^{\rho}_0 r dr d\rho  -\frac{1}{8} \int^{+\infty}_1 \rho^{1+\alpha} e^u \int^1_0 r dr d\rho \\
=&\;\frac{1}{4}\int^1_0 \rho^{3+\alpha} \ln \rho e^u d\rho - \frac{1}{16}\int^1_0 \rho^{5+\alpha} e^u d\rho -\frac{1}{16} \int^{\infty}_1 \rho^{1+\alpha}e^u d\rho.
\end{aligned}
$$
Finally, this yields
$$
\begin{aligned}
CT=\frac{1}{4}\int^1_0 \rho^{3+\alpha}\ln \rho e^u d\rho +\frac{1}{4} \int^{\infty}_1 r^{3+\alpha} \ln r e^u dr
&=\frac{1}{4}\int^{\infty}_0 r^{3+\alpha} \ln r e^u dr\\
&=\frac{1}{4\omega_3} \int_{\mathbb{R}^4} |x|^{\alpha} \ln |x| e^u dx.
\end{aligned}
$$
Therefore
\begin{align}
\label{1.15}
u(r) =-\frac{1}{4 \omega_3}\int_{\mathbb{R}^4} |x|^{\alpha} e^u dx \ln |x| + \frac{1}{4\omega_3} \int_{\mathbb{R}^4} |x|^{\alpha} \ln |x| e^u dx+ o(1) \quad \text{as}~~ r\to +\infty.
\end{align}
and
\begin{align*}
\lim_{s\to+\infty} w'(s)=(4+\alpha)-c_0<0, \quad \mbox{where } \; c_0=\frac{1}{4 \omega_3}\int_{\mathbb{R}^4} |x|^{\alpha} e^u dx.
\end{align*}

Moreover, consider the following energy function associated to \eqref{1.4}
$$
E(s)=w'''w'-\frac{w''^2}{2}-2w'-e^w,
$$
then $E'(s)=w'(w^{(4)}-4w''-e^w)\equiv 0.$ As $\lim_{s\to\pm \infty} w''(s)=\lim_{s\to \pm \infty}w'''(s)=0$ and $\lim_{s\to \pm \infty} w=-\infty$, there holds $(4+\alpha)^2=(4+\alpha-c_0)^2$, i.e. $c_0=8+2\alpha.$ So we get $(ii)$ of Theorem \ref{thm1.1}.

\subsection{\it Case 3: $N\geq 5$} We claim first $u(r)=O(\ln r)$ as $r \to +\infty.$ In fact, as $N \geq 5$,
$$\int^s_0 \frac{\sigma ^{N-1}}{1+\sigma^4} d\sigma \leq C (1+s^{N-4}),\quad \forall\; s > 0.$$
By \eqref{2.1}, there holds
\begin{align}
\label{new31}
|\Delta u(r)|\leq C \int^{+\infty}_r  s^{1-N}(1+s^{N-4}) ds\leq C(r^{2-N} +r^{-2}) \leq \frac{C'}{1+r^2}, \quad \forall\; r \geq 1.
\end{align}
Hence
\begin{align*}
|u(\sigma)| = \left|\int^r_0 s^{1-N} \int^s_0 \sigma^{N-1} \Delta u(\sigma) d\sigma ds\right| &\leq C+ C\int^r_0 s^{1-N} \int ^s_0 \frac{\sigma^{N-1}}{1+\sigma^2} d\sigma ds\\
&\leq C+ C\int^r_1 s^{1-N} (1+s^{N-2}) ds\\
&\leq C(1+\ln r).
\end{align*}

Using again \eqref{2.1}, $\Delta^2 u = r^\alpha e^u = O(r^{-4})$ as $r\to\infty$. Combining with \eqref{new31}, standard estimate shows that $u'''(r) = O(r^{-3})$ as $r\to\infty$. Consequently, $u'(r) = O(r^{-1})$ at infinity, so $u^{(k)}(r)=O(r^{-k})$ as $r \to \infty$ for $1\leq k \leq 3$. Let $w(s) = u(r)+ (4+\alpha) \ln r -\lambda_0$ with $s=\ln r,$ then
\begin{align}
\label{wcontrol}
\mbox{$w(s)\leq C$ by \eqref{2.1} and $w^{(k)}=O(1)$ at $\infty$ for $k=1,2,3$.}
\end{align}
Furthermore, $w$ satisfies
\begin{equation}
\label{1.17}
w^{(4)}+2(N-4) w^{(3)} +(N^2-10N+20) w'' -2(N-2)(N-4) w' = e^{\lambda_0}(e^w-1)
\end{equation}
where $\l_0 = \ln[2(4+\alpha)(N-2)(N-4)]$. Look at the associated energy function $$E(s)=w'''w'-\frac{w''^2}{2}+2(N-4)w''w'+\frac{1}{2} (N^2-10N+20)w'^2-e^{\lambda_0}\big(e^w-w-1\big).$$
We check readily that
$$
E'(s)=2(N-4) w''^2(s)+2(N-2)(N-4) w'^2\geq 0.
$$
So $\lim_{s \to \infty} E(s)= \ell \in \mathbb{R}$ exists by \eqref{wcontrol}.
To obtain more precise behavior of $w$ as $s\to \infty$, we consider two subcases.

\medskip\noindent
{\it Subcase I}. Suppose that $w$ is monotone near $\infty$. Then $\lim_{s\to\infty} w(s)=\gamma_1 \in\R$ exists. We claim then $\gamma_1 =0$. Otherwise, it follows from \eqref{1.17} that
\begin{equation}
\label{1.18}
\lim_{s\to \infty}\Big[w^{(4)}+2(N-4) w^{(3)} +(N^2-10N+20) w'' -2(N-2)(N-4) w'\Big]=\theta \neq 0.
\end{equation}
As $w=O(1)$ at $\infty$, there holds
$$
w^{(3)}+2(N-4)w''+(N^2-10N+20) w' \sim \theta s,\quad \text{as}~s\to \infty.
$$
By successive integrations, we would have $w(s)\sim \frac{\theta s^4}{24}$ as $s\to \infty$, which is impossible. So $\gamma_1 =0$.

\medskip\noindent
{\it Subcase II}. Suppose that $w'$ changes the sign infinitely many times near $\infty$. By the expression of $E'$, there holds $w', w''\in L^2(\R_+)$. Using \eqref{1.17} and \eqref{wcontrol}, we have
\begin{align*}
\int^s_0 w'''^2d\sigma = &\; \Big[w'''w''\Big]^s_0-\int^s_0w'' w^{(4)} d\sigma\\
=&\; O(1)- e^{\lambda_0} \int^s_0 (e^w-1)w'' d\sigma  + \int_0^s (N^2-10N+20)w''^2 d\sigma\\
& \; +\Big[(N-4)w''^2 -(N-2)(N-4)w'^2\Big]_0 ^s\\
= & \; O(1) -  e^{\lambda_0}\Big[(e^w-1)w'\Big]^s_0 +  e^{\lambda_0} \int^s_0 e^w w'^2 d\sigma =O(1).
\end{align*}
So $w'''\in L^2(\R_+)$. Similarly $w^{(4)}\in L^2(\R_+),$ so is $e^w-1$. Finally, there holds
$$\sum^4_{k=0} \int^{\infty}_0 \big|w^{(k)}\big|^2 d\sigma<\infty.$$
Then there exists a sequence $(s_j) \to \infty$ such that $\sum_{0\leq k \leq 4}|w^{(k)}|^2(s_j)\to 0$, hence $E(s_j)\to 0$, i.e.~$\lim_{s\to\infty} E(s) =0$. Considering any sequence $(t_j) \to \infty$ such that $w'(t_j)=0$. As
$$
E(t_j)=-\frac{w''^2(t_j)}{2}-e^{\lambda_0}[e^w-w-1](t_j).
$$
We must have $\lim_{j \to\infty}w(t_j)=0$, this means that $\lim_{s\to\infty} w(s)=0.$

\medskip
Therefore, we always have $\lim_{s\to\infty} w(s)=0.$
The proof of Theorem \ref{thm1.1} is completed. \qed

\section{Stability of radial solutions to \eqref{1.1-}}
\reset

The following Hardy inequalities will be important in our study of stability and stability at infinity of solutions $u_\beta$. See for instance \cite{M}.
$$
\int_{\R^N} |\Delta \phi|^2 dx \geq \frac{N^2(N-4)^2}{16} \int _{\R^N} \frac{\phi^2}{|x|^4} dx,\quad\forall ~\phi\in C^2_c(\R^N)\quad N\geq 5.
$$
and
$$
\int_{\R^4\setminus {\overline{B}_1}} |\Delta \phi|^2 dx \geq \int _{\R^4\setminus {\overline{B}_1}} \frac{\phi^2}{4|x|^4\ln ^2|x|} dx,\quad\forall ~\phi\in C^2_c(\R^4\setminus {\overline{B}_1}),
$$
$$
\int_{\R^3\setminus {\overline{B}_1}} |\Delta \phi|^2 dx \geq \frac{9}{16} \int _{\R^3\setminus \{0\}} \frac{\phi^2}{|x|^4} dx,\quad\forall ~\phi\in C^2_c(\R^3\setminus {\overline{B}_1}).
$$

\medskip

The nonexistence of stable solutions (not only radial table solutions) in dimensions less than $4$ can be referred to \cite{W}, so we omit the details.
\begin{lem}
Let $N\leq 4$, then the equation \eqref{1.1-} admits no stable solutions.
\end{lem}

For any $N\geq 3$, it's easy to see the stability at infinity for $u_\beta$ if $\beta<\beta_0$. Indeed, by \eqref{1.3-}, we have $u_{\beta}\leq -\frac{\beta_0-\beta}{2N}r^2$ for $\beta<\beta_0$, there holds then
$$|x|^\alpha e^u\leq r^\alpha e^{-\frac{\beta_0-\beta}{2N}r^2} \leq w_N(r)\quad \mbox{for $r$ large enough},$$
which implies the stability of $u_\beta$ $u_\beta$ in $\R^N\setminus B_{\widehat R_\beta}$ for $\widehat R_\beta$ large enough. Here $w_N(r)$ are the corresponding weight in the above Hardy inequalities, that is
 $$
 w_N(r) =
\begin{cases}
\begin{aligned}
&\frac{1}{4 r^4 \ln^2 r},&N=4,\\
&\frac{N^2(N-4)^2}{16r^4}, &N\ne 4.
\end{aligned}
\end{cases}
$$
Using Theorem \ref{thm1.1}, $u_{\beta_0} \sim -Cr$ when $N = 3$ and $u_{\beta_0}\sim -(8+2\alpha) \ln r$ when $N=4$, we conclude similarly that $u_{\beta_0}$ is stable outside suitable compact set, for $N = 3$ and $4$.

\medskip
Now we consider $u_{\beta_0}$ for $N \geq 5$. Let $5\leq N< N_{\alpha}$, we claim that
\begin{itemize}
\item[$(i)$] $u_{\beta_0}$ is unstable outside every compact set.
\item[$(ii)$] There exist $\beta_1<\beta_0$ depending only on $N, ~\alpha$ such that $u_\beta$ is stable and iff $\beta\leq\beta_1$.
\end{itemize}
By Theorem \ref{thm1.1}, there holds
\begin{align}
\label{new41}
\lim_{r\to \infty} u_{\beta_0}(r)+(4+\alpha)\ln r = \ln [2(4+\alpha)(N-2)(N-4)] =:\l_0.
\end{align}
So for every $m>0$, there exists $R_m>0$ such that
$$
u_{\beta_0}(r)\geq -(4+\alpha)\ln r+ \l_0-\frac{1}{m}, \quad \forall~r\geq R_m.
$$
Suppose that $u_{\beta_0}$ is stable outside a compact set $\mathcal{K} \subset B_{R_m}$, then for any $\phi \in C^2_c(\R^N\setminus{B_{R_m}})$,
\begin{align}
\label{new42}
\begin{split}
&\int_{\R^N\setminus{B_{R_m}}} |\Delta \phi|^2 dx - e^{-\frac{1}{m}+\l_0}\int_{\R^N\setminus{B_{R_m}}}\frac{\phi^2}{|x|^4} dx\\ \geq &\; \int_{\R^N\setminus{B_{R_m}}} |\Delta \phi|^2 dx - \int_{\R^N\setminus{B_{R_m}}}|x|^{\alpha}e^{u_{\beta_0}}\phi^2 dx\geq 0.
\end{split}
\end{align}
By the definition of $N_\alpha$, there holds $\frac{N^2(N-4)^2}{16} < e^{\l_0}$ for $5\leq N<N_{\alpha}$. We can choose $m$ large enough such that $\frac{N^2(N-4)^2}{16} < e^{-\frac{1}{m}+\l_0}$, but \eqref{new42} yields a contradiction with the optimality of the Hardy inequality, so $(i)$ is proved.

\medskip
For the point $(ii)$, using again \eqref{1.3-}, $u_\beta(r)\leq \frac{\beta-\beta_0}{2N} r^2$ for all $r>0$ and $\beta<\beta_0$. Let
$$
\Phi (r)=\frac{2N}{r^2}\ln \frac{N^2(N-4)^2}{16r^{4+\alpha}} + \beta_0.
$$
Obviously $\min_{r > 0} \Phi(r) = \beta'$ exists and $\beta' < \beta_0$. Moreover,
 $$
 r^{\alpha} e^{\frac{\beta -\beta_0}{2N} r^2}\leq \frac{N^2(N-4)^2}{16r^4} \;\;\mbox{in } \R^N \quad\iff \quad \beta \leq \beta'.
 $$
Hence $u_\beta$ is stable for $\beta\leq \beta'$ thanks to the Hardy inequality. Define $\beta_1 = \sup\Lambda$ where $$\Lambda := \left\{\beta < \beta_0| \; u_\beta ~\mbox{is stable in } \R^N\right\}.$$
Using the continuity and the monotonicity of $u_\beta$ with respect to $\beta$, $\Lambda = (-\infty, \beta_1]$. By $(i)$, $\beta_1 < \beta_0$.

\medskip
Finally we claim that
\begin{align}
\label{new44}
u_{\beta_0}(r) \leq -(4+\alpha)\ln r + \l_0 \;\; \mbox{in } \R^N, \quad \mbox{if }\; N\geq N_{\alpha}.
\end{align}

Let $w(s) = u(e^s) + (4+\alpha)s - \l_0$ and $P_N(\nu)=\nu (\nu-2)(\nu+N-2)(\nu+N-4)$. The left hand side of \eqref{1.17} is just $P_N(\partial_s)w$, we have then
$$
P_N(\partial_s)w(s) = e^{\lambda_0} (e^w-1)\geq e^{\lambda_0} w(s),
$$
so $Q_N(\partial_s)w(s)\geq 0$ in $\R$ with $Q_N(\nu):=P_N(\nu)-e^{\lambda_0}$.

\medskip
Remark that $P_N$ is symmetric about $\nu^*=-\frac{N-4}{2}$, i.e.~$P_N(\nu) = P_N(2\nu^* - \nu)$ in $\R$, so is $Q_N$. Moreover, if $N \geq N_\alpha$, $$\lim_{\nu \to \pm\infty} Q_N(\nu)= \infty,\quad Q_N(0) =-e^{\lambda_0}<0,\quad Q_N(\nu^*)=\frac{N^2(N-4)^2}{16}-e^{\lambda_0}\geq0.$$
Hence $Q_N(\nu) = 0$ admits four real roots $\nu_1>0>\nu_2\geq \nu_3>\nu_4$, and
\begin{equation}\label{2.19}
(\partial_s-\nu_1)(\partial_s-\nu_2)(\partial_s-\nu_3)(\partial_s-\nu_4) w(s)\geq 0 \quad \mbox{in } \R.
\end{equation}

By the definition of $w$, we have
$$\lim_{s\to-\infty}[w(s)-(4+\alpha)s]=-e^{\lambda_0},~~~\lim_{s\to-\infty} w'(s)=4+\alpha,~~\lim_{s\to-\infty}w''(s)=\lim_{s\to-\infty}w'''(s) =0.$$
We get $\lim_{s\to-\infty} {e^{-\nu_is}} w^{(k)}(s)=0$ for $2 \leq i \leq 4$, $0 \leq k \leq 3$.

\medskip
Multiplying $e^{-\nu_4 s}$ and integrating \eqref{2.19} over $(-\infty,~s)$, as $\p_s(e^{-\nu_4 s}w) = e^{-\nu_4 s}(\p_s - \nu_4)w$, there holds
$$
(\partial_s-\nu_1)(\partial_s-\nu_2)(\partial_s-\nu_3) w(s)\geq 0\quad \mbox{in } \R.
$$
Repeating this device another two times over $(-\infty,~s)$, we get
$$
(\partial_s-\nu_1) w(s)\geq 0 \quad \mbox{in } \R.$$
so $e^{-\nu_1 s}w$ is increasing in $\R$. As $\lim_{s\to +\infty} e^{-\nu_1 s} w(s)=0$, we get $w(s)\leq 0$, i.e. \eqref{new44} holds true. In other worlds, if $N \geq N_\alpha$, for any $\beta\leq \beta_0$,
$$|x|^{\alpha} e^{u_\beta(x)} \leq |x|^{\alpha} e^{u_{\beta_0}(x)}\leq \frac{N^2(N-4)^2}{16|x|^{4}}\quad \mbox{in } \R^N,$$
hence $u_{\beta}$ is stable by the Hardy inequality. The proof of Theorem \ref{thm1.2} is completed.\qed

\medskip
\noindent

\end{document}